\documentclass[12pt]{article}
\usepackage{ytableau}
\usepackage{amssymb}
\usepackage[margin=1in]{geometry}
\usepackage{mathrsfs}
\usepackage{stmaryrd}
\usepackage{amsfonts,amsmath,amssymb,amscd}
\usepackage{shadow}
\usepackage{graphicx}
\usepackage{color}
\usepackage{longtable}
\usepackage{pstricks,multido}
\usepackage{hyperref}
\usepackage{qtree}
\allowdisplaybreaks

\newtheorem{thm}{Theorem}[section]
\newtheorem{lem}[thm]{Lemma}

\newtheorem{cor}[thm]{Corollary}

\newtheorem{conj}{Conjecture}[section]

\def\qed{\hfill \rule{4pt}{7pt}}

\def\pf{\noindent {\it{Proof.} \hskip 2pt}}

\numberwithin{equation}{section}

\begin{document}
\begin{center}
{\large\bf  On the Enumeration of
$(s,s+1,s+2)$-Core Partitions}
\end{center}

\begin{center}

{\small Jane Y.X. Yang$^1$, Michael X.X. Zhong$^2$, Robin D.P. Zhou$^3$}

Center for Combinatorics, LPMC-TJKLC\\
Nankai University\\
Tianjin 300071, P.R. China

$^1$yangyingxue@mail.nankai.edu.cn,
$^2$michaelzhong@mail.nankai.edu.cn,
$^3$robin@cfc.nankai.edu.cn

\end{center}

\begin{abstract}
Anderson established a connection between core partitions and order ideals of certain posets by mapping a partition to its $\beta$-set.
In this paper, we give a characterization of the poset $P_{(s,s+1,s+2)}$ whose order ideals correspond to   $(s,s+1,s+2)$-core partitions.
Using this characterization, we obtain
 the number of $(s,s+1,s+2)$-core partitions, the maximum size and the average size of an $(s,s+1,s+2)$-core partition,
confirming three conjectures posed by Amdeberhan.
\end{abstract}

\noindent {\bf Keywords}: core partition, hook length, $\beta$-set, poset, order ideal

\noindent {\bf AMS  Subject Classifications}:  05A15, 05A17, 06A07

\section{Introduction}

The objective of this paper is to prove three conjectures of Amdeberhan on
$(s,s+1,s+2)$-core partitions.

A partition $\lambda$ of a positive integer $n$ is a finite nonincreasing sequence of positive integers
$(\lambda_1, \lambda_2, \ldots, \lambda_m)$ such that
$\lambda_1+\lambda_2+\cdots+\lambda_m=n$.
We write
$\lambda=(\lambda_1,\lambda_2,\ldots,\lambda_m)\vdash n$
 and we say that $n$ is the size of $\lambda$ and $m$ is the
 length of $\lambda$.
The Young diagram of $\lambda$ is defined to be a left-justified array
of $n$ boxes with $\lambda_i$ boxes in the $i$th row.
For each box $B$ in $\lambda$, it determines a hook
consisting of the box $B$ itself and
 boxes directly to the right and directly below $B$.
 The hook length of $B$, denoted $h(B)$, is the number
 of boxes in the hook of $B$.

For a partition $\lambda$,   the $\beta$-set of $\lambda$, denoted $\beta(\lambda)$, is defined to be the set of hook lengths of the boxes in the first column of $\lambda$.
For example, Figure
\ref{1.1fig} illustrates the Young diagram and the  hook lengths of a partition $\lambda=(5,3,2,2,1)$. The
$\beta$-set of $\lambda$ is $\beta(\lambda)=\{9,6,4,3,1\}$.
Notice that a partition $\lambda$ is uniquely determined by its $\beta$-set.
Given a decreasing sequence of positive integers $(h_1,h_2,\ldots, h_m)$,
it is easily seen that the unique partition $\lambda$ with $\beta(\lambda)=\{h_1,h_2,\ldots, h_m\}$ is
\begin{equation}\label{betaset}
\lambda=(h_1-(m-1), h_2-(m-2), \dots, h_{m-1}-1, h_m).
\end{equation}

\begin{figure}[h]
\begin{center}
 \begin{ytableau}
     9& 7 & 4 & 2&1\\
     6 & 4 & 1\\
    4 & 2 \\
     3 & 1\\
     1
 \end{ytableau}
 \end{center}
\caption{The Young diagram of $\lambda=(5,3,2,2,1)$.}
\label{1.1fig}
\end{figure}

For a positive integer $t$, a partition $\lambda$ is a $t$-core partition, or simply a $t$-core,  if
it contains no box whose hook length is a multiple of
 $t$.
Let $s$ be a positive integer distinct with $t$, we say that
$\lambda$ is an $(s,t)$-core if
it is simultaneously an $s$-core and a $t$-core.
For example, the partition $\lambda=(5,3,2,2,1)$ in
Figure \ref{1.1fig}  is a $(5,8)$-core.
In general, an $(a_1,a_2,\ldots, a_r)$-core partition can be defined for distinct positive integers $a_1,a_2,\ldots, a_r$.
Since a $t$-core is an $s$-core if $s$ is a multiple of $t$, we assume that there is no element in $\{a_1,a_2,\ldots, a_r\}$ that is a multiple of another element.

Let $s$ and $t$ be two coprime positive integers.
Anderson \cite{Anderson} showed that the number of $(s,t)$-core partitions equals
${s+t\choose s}/(s+t)$.
 Ford, Mai and Sze \cite{Ford}  proved that the number of self-conjugate $(s,t)$-core partitions equals ${\lfloor \frac{s}{2}\rfloor+\lfloor\frac{t}{2}\rfloor \choose \lfloor\frac{s}{2}\rfloor}$.
Furthermore,   Olsson and  Stanton \cite{Olsson} proved
that there exists a unique $(s,t)$-core partition with the maximum size  $(s^2-1)(t^2-1)/24$. A simpler proof was
 provided by Tripathi \cite{Tripathi}.
Armstrong, Hanusa and Jones \cite{Armstrong}  conjectured that the average size of an $(s,t)$-core partition and the average size of a self-conjugate $(s,t)$-core are both equal
$(s+t+1)(s-1)(t-1)/24$.
Stanley and Zanello \cite{Stanley} showed that the average size of
 an $(s,s+1)$-core equals ${s+1\choose 3}/2$. Chen, Huang and Wang \cite{Chen} proved the conjecture for the average size of a self-conjugate $(s,t)$-core.

Concerning the enumeration of $(s,s+1,s+2)$-core partitions, Amdeberhan \cite{Amdeberhan} posed  three conjectures.

\begin{conj} \label{conj1}
Let $C_k$ be the $k$th Catalan number, that is, $C_k=\frac{1}{k+1}{2k\choose k}$. Let $s$ be a
positive integer, the number $r(s)$ of $(s,s+1,s+2)$-core partitions equals
\begin{equation*}
\sum_{k\geq 0}{s \choose 2k}C_k.
\end{equation*}
\end{conj}

\begin{conj}\label{conj2}
Let $s$ be a positive integer, the size $l(s)$ of the largest $(s,s+1,s+2)$-core partition equals   \begin{equation*}
l(s)=\left\{
       \begin{array}{ll}
         m{m+1\choose 3}, & \mbox{if }  s=2m-1, \\[6pt]
         (m+1){m+1\choose 3}+{m+2\choose 3},  & \mbox{if } s=2m.
       \end{array}
     \right.
\end{equation*}
\end{conj}

\begin{conj}\label{conj3}
Let $s$ be a positive integer, the sum $h(s)$ of the sizes of all   $(s,s+1,s+2)$-core partitions equals
\begin{equation*}
h(s)=\sum _{j=0}^{s-2} {j+3\choose 3}\sum _{i=0}^{  \lfloor j/2 \rfloor}
{j\choose 2\,i}C_i.
\end{equation*}
Equivalently, the average size of an $(s,s+1,s+2)$-core partition is $\frac{h(s)}{r(s)}$.
\end{conj}

Anderson \cite{Anderson} characterized the $\beta$-sets of $(s,t)$-core partitions as order ideals of a poset $P_{(s,t)}$,
where
\[P_{(s,t)}=\mathbb{N}^+\setminus \{n \in\mathbb{N}^+\mid n=k_1s+k_2t \mbox{ for some } k_1,k_2\in \mathbb{N}\}\]
and $y\geq x$ in $P_{(s,t)}$ if there exist $y=y_0,y_1,y_2,\ldots,y_l=x\in P_{(s,t)}$ such that $y_i-y_{i+1}\in\{s,t\}$.
We show that the above characterization can be
generalized to $(a_1, a_2, \ldots, a_r)$-core partitions.
More precisely, for   positive integers $a_1,a_2,\ldots,a_r$, we define
\[P_{(a_1,a_2,\ldots, a_r)}=\mathbb{N}^+\setminus \{n \in\mathbb{N}^+\mid n=k_1a_1+k_2a_2+\cdots+k_ra_r \mbox{ for some } k_1,k_2,\ldots,k_r\in \mathbb{N}\},\]
where
$y\geq x$ in $P_{(a_1,a_2,\ldots, a_r)}$ if there exist $y=y_0,y_1,y_2,\ldots,y_l=x\in P_{(a_1,a_2,\ldots, a_r)}$ such that $y_i-y_{i+1}\in\{a_1,a_2,\ldots, a_r\}$.
It can be shown  that   $\beta$-sets of $(a_1,a_2,\ldots, a_r)$-core partitions are exactly   order ideals of the poset $P_{(a_1,a_2,\ldots, a_r)}$.
Based on this characterization, we shall prove the above three conjectures.

\section{Proof of Conjecture \ref{conj1}}\label{sec:2}

In this section, we show that a partition is an $(a_1, a_2, \ldots, a_r)$-core if and only if its $\beta$-set is an order ideal of the poset $P_{(a_1, a_2, \ldots, a_r)}$.
We shall use this correspondence to derive a formula for the number of $(s,s+1,s+2)$-core partitions.


Let $P$ be a poset. For two elements $x,y$ in $P$, we say $y$ covers $x$ if $x< y$ and there exists no element $z\in P$
satisfying $x< z< y$.
 The Hasse diagram of a finite poset $P$ is a graph whose vertices are the elements of $P$, whose edges are the cover relations, and such that if $y$ covers $x$ then there is an edge connecting $x$ and $y$ and
 $y$ is placed above $x$. An order ideal of $P$ is
a subset $I$ such that  if any $y\in I$ and $x\leq y$ in $P$, then $x\in I$. Let $J(P)$ denote the set of order ideals of  $P$, see Stanley \cite{Stanleybook}.

In the following theorem, Anderson \cite{Anderson} established a correspondence between core partitions and
 order ideals of a certain poset by mapping a partition to its $\beta$-set.

\begin{thm}\label{thm:core-poset}
Let $s,t$ be two coprime positive integers, and let $\lambda$ be a partition of $n$. Then $\lambda$ is an $s$-core (or $(s,t)$-core) partition if and only if $\beta(\lambda)$ is an order ideal of $P_s$ (or $P_{(s,t)}$).
\end{thm}

For example, let $s=3$ and $t=4$.
We can construct all $(3,4)$-core partitions by finding order ideals of $P_{(3,4)}$.
It is easily checked that $P_{(3,4)}=\{1,2,5\}$ with the
partial order $5>2$ and $5>1$. Hence the order ideals of $P_{(3,4)}$ are
$\emptyset$, $\{1\}$, $\{2\}$, $\{2,1\}$ and $\{5,2,1\}$.
The corresponding $(3,4)$-core partitions
are $\emptyset$, $(1)$, $(2)$, $(1,1)$ and $(3,1,1)$, respectively.

 Theorem \ref{thm:core-poset} can be extended to $(s,s+1,s+2)$-core partitions.

\begin{thm}\label{thm:general}
Let $a_1,a_2,\ldots, a_r$ be a sequence of positive integers, and let $\lambda$ be a partition of $n$. Then $\lambda$ is an $(a_1,a_2,\ldots, a_r)$-core if and only if $\beta(\lambda)$ is an order ideal of $P_{(a_1,a_2,\ldots, a_r)}$.
\end{thm}

\pf
Assume that $\lambda$ is an $(a_1,a_2,\ldots, a_r)$-core, we proceed to prove that $\beta(\lambda)$ is an order ideal of $P_{(a_1,a_2,\ldots, a_r)}$.
First, we claim that $\beta(\lambda)$ is a subset of
$P_{(a_1,a_2,\ldots, a_r)}$. Otherwise, suppose that $h$ is an element in $\beta(\lambda)$ but it is not contained in
$P_{(a_1,a_2,\ldots, a_r)}$.
 By the definition of $P_{(a_1,a_2,\ldots, a_r)}$, there exist nonnegative integers $k_1,k_2,\ldots,k_r$  such that
\[h=k_1a_1+k_2a_2+\cdots+k_ra_r.\]
Without loss of generality, we may assume that $k_1>0$.
Since $\lambda$ is an $(a_1,a_2,\ldots, a_r)$-core partition, it is an $a_r$-core partition. By
Theorem \ref{thm:core-poset}, we see that $\beta(\lambda)$ is an order ideal of $P_{a_r}$.
Since $k_1a_1+k_2a_2+\cdots+k_{r-1}a_{r-1}\in P_{a_r}$, it is easily seen that $k_1a_1+k_2a_2+\cdots+k_{r-1}a_{r-1}\in \beta(\lambda)$.
Now, since $\lambda$ is an $a_{r-1}$-core partition,
 we find that $k_1a_1+k_2a_2+\cdots+k_{r-2}a_{r-2}\in \beta(\lambda)$.
Continuing the above process, we eventually obtain that $k_1a_1\in \beta(\lambda)$, contradicting the fact that $\lambda$ is an $a_1$-core partition.
Thus the claim is proved.

To prove that $\beta(\lambda)$ is an order ideal of $P_{(a_1,a_2,\ldots, a_r)}$, we assume that
 $y\in \beta(\lambda)$ and $x$ is covered by $y$ in
$P_{(a_1,a_2,\ldots, a_r)}$. We need  to show that $x\in \beta(\lambda)$.
Since $y$ covers $x$ in $P_{(a_1,a_2,\ldots, a_r)}$, there exists  $1\leq i\leq r$ such that $y-x=a_i$.
From the fact that $\beta(\lambda)$ is an order ideal of $P_{a_i}$, we see that $x\in \beta(\lambda)$.

Conversely, assume that $\lambda$ is a partition such that $\beta(\lambda)$ is an order ideal of  $P_{(a_1,a_2,\ldots, a_r)}$. We aim to show
that $\lambda$ is an $(a_1,a_2,\ldots, a_r)$-core partition.
We now claim that $\lambda$ is an $a_1$-core partition.
By Theorem \ref{thm:core-poset}, it suffices to prove
that $\beta(\lambda)$ is an order ideal of $P_{a_1}$.
Notice that $\beta(\lambda)$ is  a subset of $P_{a_1}$ since
$P_{(a_1,a_2,\ldots, a_r)}\subseteq P_{a_1}$.
To prove that $\beta(\lambda)$ is an order ideal of $P_{a_1}$, we assume that $y\in \beta(\lambda)$, $x\in P_{a_1}$ and $y-x=a_1$.
It remains to show that $x\in\beta(\lambda)$.
First, we show that $x\in P_{(a_1,a_2,\ldots, a_r)}$.
Otherwise, we assume that there exist nonnegative integers $c_1,c_2,\ldots,c_r$ such that
\[x=y-a_1=c_1a_1+c_2a_2+\cdots+c_ra_r,\]
or equivalently,
\[y=(c_1+1)a_1+c_2a_2+\cdots+c_ra_r.\]
It follows that $y\not \in P_{(a_1,a_2,\ldots, a_r)}$,
which contradicts the assumption $y\in P_{(a_1,a_2,\ldots, a_r)}$.
 So we have $x\in P_{(a_1,a_2,\ldots, a_r)}$.
Since $\beta(\lambda)$ is an order ideal of $P_{(a_1,a_2,\ldots, a_r)}$ and $y-x=a_1$, we obtain $x\in \beta(\lambda)$.
Thus, $\beta(\lambda)$ is an order ideal of $P_{a_1}$, which implies that $\lambda$ is an $a_1$-core.
This proves the claim.

Similarly, it can be shown that $\lambda$ is
an $a_i$-core for $2\leq i\leq r$.
Hence  $\lambda$ is an $(a_1,a_2,\ldots, a_r)$-core.
This completes the proof.
\qed

Theorem \ref{thm:general} establishes a correspondence between
$(s,s+1,s+2)$-core partitions and order ideals of
$P_{(s,s+1,s+2)}$.
The following description of $P_{(s,s+1,s+2)}$ can be
used to compute the number of order ideals of
$P_{(s,s+1,s+2)}$.
For convenience, we denote $P_{(s,s+1,s+2)}$ by $T_s$.
Given positive integers $a\leq b$,
we denote $\{a,a+1,\ldots,b\}$ by $[a,b]$.

\begin{thm}\label{thm:poset}
Let $s\geq 3$ be a positive integer. Then
$T_s$ is graded of length $\lfloor\frac{s}{2}\rfloor-1$.
More precisely, we have
\[T_s=B_{0}\cup B_{1} \cup \cdots \cup B_{\lfloor\frac{s}{2}\rfloor-1},\]
where $B_{k}=[1+k(s+2),(k+1)s-1]$ denotes the set of the elements with rank $k$.
For $1\leq k\leq \lfloor\frac{s}{2}\rfloor-1$, each element $b$ in $B_{k}$ covers exactly the three elements $b-s,b-(s+1),b-(s+2)$ in $B_{k-1}$.
\end{thm}

\pf
By the definition of $P_{(s,s+1,s+2)}$, it is easily seen that
\[T_s=P_{(s,s+1,s+2)}=B_{0}\cup B_{1} \cup \cdots \cup B_{\lfloor\frac{s}{2}\rfloor-1}.\]
We proceed to show that $T_s$ is graded.
Examining the definition of  $T_s$,
we see that for each element $b$ in $B_{k}$, the possible elements covered by $b$ are $b-s,b-(s+1),b-(s+2)$.
Since $b\in B_k=[1+k(s+2),(k+1)s-1]$, it is easily checked that
all the elements $b-s,b-(s+1),b-(s+2)$ are  in $B_{k-1}=[1+(k-1)(s+2),ks-1]$ for $k \geq 1$.
Conversely, either $b+s$ or $b+(s+2)$ is in $B_{k+1}$
for $k< \lfloor \frac{s}{2} \rfloor-1$,
that is, $b$ must be covered by at least one element in $B_{k+1}$.
Hence $T_s$ is graded of length $\lfloor\frac{s}{2}\rfloor-1$.
This completes the proof.
\qed

\begin{figure}[h]
\centerline{\includegraphics[width=15cm]{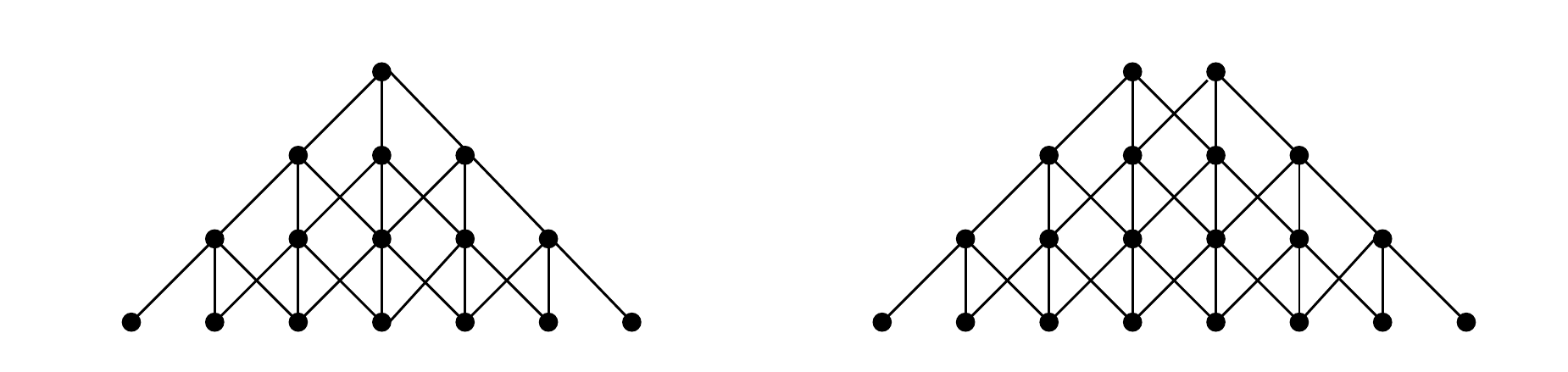}}
\caption{The Hasse diagrams of the posets $T_8$ and $T_9$.}
\label{figT8-9}
\end{figure}

According to Theorem \ref{thm:poset}, the Hasse diagram of $T_s$ can be easily constructed.
For example, Figure \ref{figT8-9} illustrates the Hasse diagrams of the posets $T_8$ and $T_9$.

  Theorem \ref{thm:poset} enables us to
compute the number
of order ideals of $T_s$.
To this end, we shall partition $J(T_s)$ according to the smallest missing element of rank $0$ in an order ideal.
Note that the elements of rank $0$ in $T_s$ are just the minimal elements.
For $1\leq i\leq s-1$, let $J_i(T_s)$ denote the set of order ideals of $T_s$ such that $i$ is the smallest missing element of rank $0$.
Let $J_s(T_s)$ denote the set of order ideals which contain all minimal elements in $T_s$.
Then we can write $J(T_s)$ as \[ J(T_s)=\bigcup_{i=1}^sJ_i(T_s).\]

\begin{figure}[h]
\centerline{\includegraphics[bb={115  552  478 743},width=15cm]{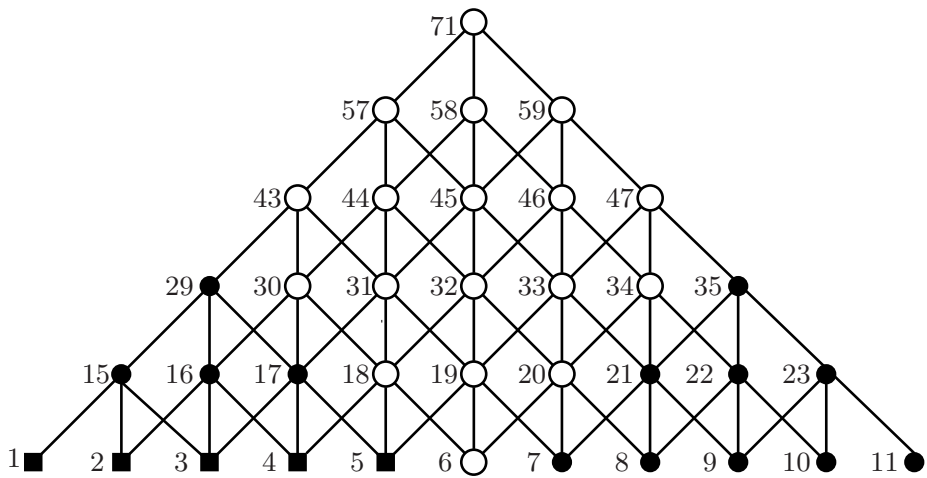}}
\caption{The  elements of an order ideal $I\in J_6(T_{12})$.}
\label{1.3fig}
\end{figure}

Figure \ref{1.3fig} gives an illustration of the   elements contained in an order ideal in $J_6(T_{12})$.
We see that an order ideal $I\in J_6(T_{12})$ must contain the elements labeled by squares, but does not contain any elements represented by open  circles.
The elements represented by solid  circles  may or may not  appear in $I$.
That is, $I$ can be decomposed into three parts, one is $\{1,2,3,4,5\}$, one is isomorphic to an order ideal of $T_4$ and one is isomorphic to an order ideal of $T_6$.

In general, for $2\leq i\leq s$ and an order ideal $I\in J_i(T_s)$, we can decompose it into three parts, one is $\{1,2,\ldots,i-1\}$, one is isomorphic to an order ideal of $T_{i-2}$ and one is isomorphic to an order ideal of $T_{s-i}$.
We shall use  this decomposition to prove Conjecture
\ref{conj1}. Recall that the
Motzkin number \cite{A001006} $M_s$ equals
\[
\sum_{k\geq 0}{s \choose 2k}C_k.
\]

By Theorem \ref{thm:general}, to prove Conjecture \ref{conj1}, it suffices to show that the number $r(s)$ of order ideals of $T_s$ equals $M_s$.

\noindent{\it Proof of Conjecture \ref{conj1}.}
It is easily checked that the conclusion is correct
when $s=0,1,2$.
Suppose now $s\geq 3$.
For an order ideal $I\in J_1(T_s)$, $I$ is isomorphic to
an order ideal of $T_{s-1}$.
For $2 \leq i\leq s$ and an order ideal $I\in J_i(T_s)$,
$I$ can be decomposed into
three parts,
one is $\{1,2,\ldots,i-1\}$, one is isomorphic to an order ideal of $T_{i-2}$ and one is isomorphic to an order ideal of $T_{s-i}$.
Hence we have
\begin{equation}\label{equ:2.1}
r(s)=r(s-1)+\sum_{i=2}^sr(i-2)r(s-i).
\end{equation}
It is known that the Motzkin number $M_s$ satisfies
 recurrence relation (\ref{equ:2.1}) with the same initial conditions as $r(s)$. This yields that $r(s)=M_s$, and hence the proof is complete.
\qed

\section{Proof of Conjecture \ref{conj2}}\label{sec:3}

In this section, we show that the partition $\kappa_s$ corresponding to the order ideal consisting of all elements in the poset $T_s$ is of maximum size.
Moreover, we show that if $s$ is even, then $\kappa_s$ is the unique $(s,s+1,s+2)$-core partition of maximum size,
and if $s\geq 3$ is odd, then there is exactly another $(s,s+1,s+2)$-core partition of maximum size which is the conjugate of $\kappa_{s}$.
This leads to a proof of Conjecture \ref{conj2}.

We need the following three lemmas to
characterize  order ideals of $T_s$ corresponding to   $(s,s+1,s+2)$-core partitions of maximum size.

Recall that for an order ideal $\beta=\{h_1,h_2,\ldots,h_m\}$ of $T_s$ where the elements are listed in decreasing order,
the corresponding $(s,s+1,s+2)$-core partition $\lambda$ is given by $\lambda=(h_1-(m-1),h_2-(m-2),\ldots,h_m)$,
whose size is given by
\begin{align}\label{equ:size}
|\lambda|=\sum_{i=1}^m h_{i}-{m \choose 2}.
\end{align}

For example,  $\beta=\{16,15,4,3,2,1\}$ is an order ideal of $T_{12}$, which corresponds to a $(12,13,14)$-core partition
$\lambda=(11,11,1,1,1,1)$ of size $26$.

Note that $B_k$ is the set of elements in $T_s$ of rank $k$,
that is, \[B_k=[1+k(s+2),(k+1)s-1].\]

\begin{lem}\label{lem:continuity}
Let $\lambda$ be an $(s,s+1,s+2)$-core partition of maximum size.
If $\beta(\lambda)$ contains an element $i$
 that is in $B_k$, then $\beta(\lambda)$ contains all the
 elements in $[i,(k+1)s-1]$.
\end{lem}

\begin{figure}[h]
\centerline{\includegraphics[width=15cm]{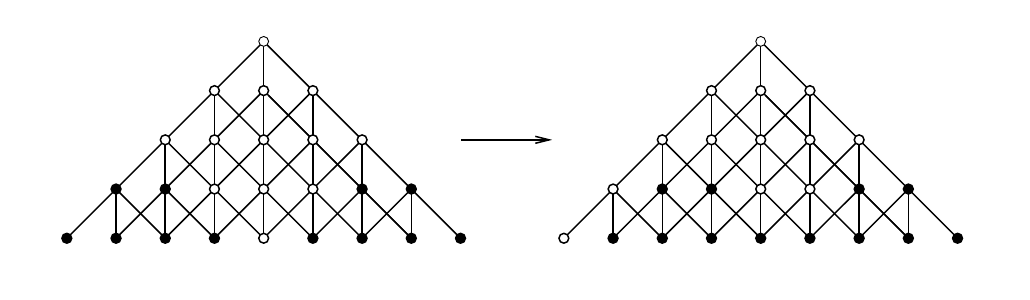}}
\caption{An example in $T_{10}$.}
\label{fig:lemma1}
\end{figure}

\pf
Assume to the contrary that the lemma is not valid, that is,
there exist elements $i,j\in B_k$ such that $i<j$, $i\in \beta(\lambda)$ and $j\not \in \beta(\lambda)$.
We choose $k$ to be the smallest and  $i$ to be the smallest after
 $k$ is chosen.
For any $p\in\beta(\lambda)$ such that $p\geq q$ in $T_s$ for some $q\in[i,j-1]$, we replace it by $p+1$, see Figure \ref{fig:lemma1} for an illustration.
This leads us to a new order ideal $\beta'$ with the same cardinality as $\beta(\lambda)$ and a larger sum of the elements.
By relation \eqref{equ:size}, the size of the $(s,s+1,s+2)$-core partition corresponding to $\beta'$ is larger than that of $\lambda$, which contradicts   the assumption that $\lambda$ is of maximum size.
This proves the lemma.
\qed

\begin{figure}[h]
\centerline{\includegraphics[width=15cm]{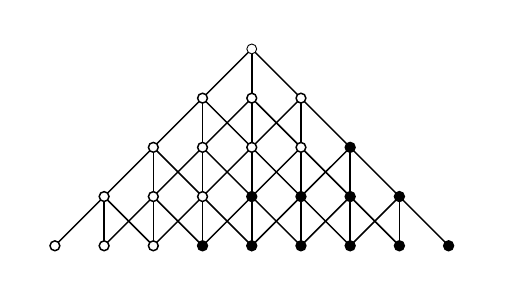}}
\caption{The order ideal $\beta_{4,2}$ of $T_{10}$.}
\label{fig:lemma2}
\end{figure}

For $1\leq i\leq s-1$, let $\beta_{i,0}$ be the unique order ideal in $T_s$ which is isomorphic to $T_{s-i}$ and contains all the elements in $[i+1,s-1]$.
For $1\leq j\leq \lfloor \frac{s-i+1}{2} \rfloor$, let $\beta_{i,j}$ be the union of $\beta_{i,0}$ and the chain consisting of $i,i+(s+2),\ldots,i+(j-1)(s+2)$.
 For example, the order ideal $\beta_{4,2}$ of $T_{10}$ is given in Figure \ref{fig:lemma2}.

For $ 1\leq i\leq s-1$, $0\leq j\leq \lfloor\frac{s-i+1}{2}\rfloor$,
$\beta_{i,j}$ is an order ideal of $T_s$.
Let $\lambda_{i,j}$ be the unique partition such that $\beta(\lambda_{i,j})=\beta_{i,j}$.
By Theorem \ref{thm:general}, for each $\beta_{i,j}$, $\lambda_{i,j}$ is an $(s,s+1,s+2)$-core partition.
Let $\lambda$ be an $(s,s+1,s+2)$-core partition of maximum size.
We shall show that $\lambda$  equals $\lambda_{i,j}$ for some integers $i,j$.
From Lemma \ref{lem:continuity}, we get that $s-1\in \beta(\lambda)$, so that there exists an integer $i$ such that $[i,s-1]$ is contained in $\beta(\lambda)$.

\begin{lem}\label{lem:triangleideal}
Assume that $s\geq 3$. Let $\lambda$ be an $(s,s+1,s+2)$-core partition of maximum size.
Then there exist some integers $1\leq i\leq s-1$ and $0\leq j\leq \lfloor\frac{s-i+1}{2} \rfloor$ such that
$\lambda=\lambda_{i,j}$.
\end{lem}

\pf
Let $i$  be the minimal integer such that $[i,s-1]$ is contained in $\beta(\lambda)$ and $j$ the maximal integer such that $i+(j-1)(s+2)\in \beta(\lambda)$.
We proceed to show that $\lambda=\lambda_{i,j}$, or equivalently, $\beta(\lambda)=\beta_{i,j}$.

By the choice of $i$ and $j$, the proof of Lemma \ref{lem:continuity} shows that  $\beta(\lambda)\subseteq \beta_{i,j}$.
Hence it remains to show that $\beta_{i,j}\subseteq\beta(\lambda)$.
Assume to the contrary that $\beta_{i,j}\not\subseteq\beta(\lambda)$,
that is, there exists an element in $\beta_{i,j}$ which is not contained in $\beta(\lambda)$.
Let $p$ be the smallest element such that $p\in \beta_{i,j}$ and
$p\not \in \beta(\lambda)$.

Let $\beta'$ denote the set $\beta(\lambda)\cup \{p\}\setminus \{i+(j-1)(s+2)\}$. We claim that $\beta'$ is an order
ideal of $T_s$ and it corresponds to an $(s,s+1,s+2)$-core partition of size larger than $|\lambda|$.

First, we show that  $\beta'$ is an order ideal of $T_s$.
Let $\gamma=\beta(\lambda)\cup \{p\}$.
To prove that $\beta'$ is an order ideal of $T_s$,
it is sufficient to show that $\gamma$ is an order ideal of
$T_s$ and $i+(j-1)(s+2)$ is a maximal element in $\gamma$.
Let $q$ be an arbitrary element of $T_s$ such that $q<p$ in the poset $T_s$.
Notice that $\beta_{i,j}$ is an order ideal of $T_s$.
 This implies that $q\in \beta_{i,j}$ since $p\in \beta_{i,j}$ and $q<p$ in $T_s$.
By the choice of $p$, we see that  $q\in \beta(\lambda)\subseteq \gamma$.
Hence $\gamma$ is an order ideal of $T_s$.

To prove that $\beta'$ is an order ideal of $T_s$,
it remains to show that $i+(j-1)(s+2)$ is a maximal element
of the order ideal $\gamma$.
By the definition of $\beta_{i,j}$,    $i+(j-1)(s+2)$ is a maximal element
of  $\beta_{i,j}$.
Since $\gamma\subseteq \beta_{i,j}$ and $i+(j-1)(s+2)\in \gamma$,
we obtain that $i+(j-1)(s+2)$ is a maximal element of $\gamma$.
Hence $\gamma$ is an order ideal of
$T_s$ and $i+(j-1)(s+2)$ is a maximal element in $\gamma$.
 So we deduce that $\beta'$ is an order ideal of $T_s$.

Let $\mu$ be the partition determined by $\beta(\mu)=\beta'$.
 By Theorem \ref{thm:general}, $\mu$ is an $(s,s+1,s+2)$-core.
We aim to show that $|\mu|>|\lambda|$.
Because of relation \eqref{equ:size},  it suffices to show that $p>i+(j-1)(s+2)$.
Assume to the contrary that $p\leq i+(j-1)(s+2)$.
Since $p\in \beta_{i,j}$, we obtain that
$p \in B_{k}\cap \beta_{i,j}=[i+k(s+2),(k+1)s-1]$ for some
integer $0\leq k \leq j-1$.
Notice that $i+(j-1)(s+2)\in \beta(\lambda)$ and
$\beta(\lambda)$ is an order ideal of $T_s$.
Since $i+k(s+2)\leq i+(j-1)(s+2)$ in $T_s$, we have $i+k(s+2)\in \beta(\lambda)$.
From Lemma \ref{lem:continuity} we see that
$[i+k(s+2),(k+1)s-1]\subseteq \beta(\lambda)$.
It follows that $p\in \beta(\lambda)$,
which contradicts the assumption that $p\not \in \beta(\lambda)$.
Thus $p>i+(j-1)(s+2)$, that is, $|\mu|>|\lambda|$,
contradicting the condition that $\lambda$ is of maximum size.
This proves that $\beta_{i,j}\subseteq \beta(\lambda)$.
 So we conclude that $\lambda=\lambda_{i,j}$, and this completes
 the proof.
\qed

\begin{lem}\label{lem:triangle}
Given $1\leq i\leq s-1$,
we have $|\lambda_{i,j}| \leq| \lambda_{i,\lfloor\frac{s-i+1}{2}\rfloor}|$ for $0\leq j\leq\lfloor\frac{s-i+1}{2}\rfloor$, with the equality holds if and only if $j=\lfloor\frac{s-i+1}{2}\rfloor$, or $s$ is odd, $i=1$ and $j=0$.
\end{lem}

\pf
By relation \eqref{equ:size}, the size of $\lambda_{i,j}$ equals
\begin{align}\label{equ:lambda(i,j)}
|\lambda_{i,j}|& = \sum_{h\in \beta_{i,j}}{h}-{|\beta_{i,j}| \choose 2}\notag\\[6pt]
& =\sum_{h\in \beta_{i,0}}{h}+\sum_{p=0}^{j-1}{\big(i+p(s+2)\big)}-{|\beta_{i,0}|+j \choose 2}\notag\\[6pt]
& =\sum_{h\in \beta_{i,0}}{h}-{|\beta_{i,0}| \choose 2}+\sum_{p=0}^{j-1}{\big(i+p(s+2)\big)}-{|\beta_{i,0}|+j \choose 2}+{|\beta_{i,0}| \choose 2}\notag\\[6pt]
& =|\lambda_{i,0}|+ij+(s+1){j \choose 2}-|\beta_{i,0}|j.
\end{align}
By the definition of $\beta_{i,0}$ and Theorem \ref{thm:poset}, we obtain that
\[
|\beta_{i,0}|=|T_{s-i}|=\left\{
\begin{array}{ll}
k^2+k, & \hbox{if } i=s-2k-1,\\[6pt]
k^2, & \hbox{if } i=s-2k.
\end{array}
\right.
\]
Hence
\begin{align}\label{lambda(i,j)}
|\lambda_{i,j}|=\left\{\begin{array}{ll}
|\lambda_{i,0}|+ij+(s+1){j \choose 2}-(k^2+k)j,& \hbox{if } i=s-2k-1,\\[6pt]
|\lambda_{i,0}|+ij+(s+1){j \choose 2}-k^2j, & \hbox{if }  i=s-2k.
\end{array}
\right.
\end{align}
In particular, for $j=\lfloor\frac{s-i+1}{2}\rfloor$, we have
\begin{align}\label{lambda(i,j,0)}
|\lambda_{i,\lfloor\frac{s-i+1}{2}\rfloor}|-|\lambda_{i,0}|=\left\{\begin{array}{ll}
(\frac{ik}{2}+i)(k+1),& \hbox{if } i=s-2k-1,\\[6pt]
\frac{(i-1)(k^2+k)}{2}, & \hbox{if }  i=s-2k,
\end{array}
\right.
\end{align}
which implies that
\begin{equation}\label{equ:lambda(i,j)>lambda(i,0)}
|\lambda_{i,\lfloor\frac{s-i+1}{2}\rfloor}| \geq |\lambda_{i,0}|.
\end{equation}

For fixed integers $i$ and $s$,
we see that $|\lambda_{i,j}|$ is a quadratic function of $j$ with a positive leading coefficient.
Hence the maximum value of $|\lambda_{i,j}|$ is obtained at $j=0$ or $j=\lfloor\frac{s-i+1}{2}\rfloor$ when $j$  ranges over
$[0,\lfloor\frac{s-i+1}{2}\rfloor]$.
In view of  \eqref{equ:lambda(i,j)>lambda(i,0)},
we conclude  that
\begin{equation}\label{inequa1}
|\lambda_{i,j}| \leq |\lambda_{i,\lfloor\frac{s-i+1}{2}\rfloor}|
\end{equation}
for $0\leq j\leq \lfloor\frac{s-i+1}{2}\rfloor$.
Moreover, we have
\begin{equation*}
|\lambda_{i,j}| < |\lambda_{i,\lfloor\frac{s-i+1}{2}\rfloor}|
\end{equation*}
for $0< j< \lfloor\frac{s-i+1}{2}\rfloor$.
Hence  \eqref{inequa1}   holds with equality
only when $j=0$
or $j=\lfloor\frac{s-i+1}{2}\rfloor$.
Assume that \eqref{inequa1}  holds with equality for $j=0$,
that is $|\lambda_{i,0}| = |\lambda_{i,\lfloor\frac{s-i+1}{2}\rfloor}|$.
It follows from  \eqref{lambda(i,j,0)} that $i=1$ and $s=2k+1$ for some $k$. Conversely, for $i=1$, $j=0$ and $s=2k+1$, by  \eqref{lambda(i,j,0)} we see that \eqref{inequa1} holds with equality, namely, $|\lambda_{1,0}| = |\lambda_{1,\lfloor\frac{s}{2}\rfloor}|$.
This completes the proof.
\qed

The following theorem provides a
characterization of  $(s,s+1,s+2)$-core partitions of maximum size
in terms of corresponding order ideals of $T_s$ under the map
$\beta$. We shall use the common notation $\lambda'$ for
the conjugate of a partition $\lambda$.

\begin{thm}\label{thm:maxcore}
Assume that $s\geq 3$. Let $\kappa_s$ be the $(s,s+1,s+2)$-core partition
such that $\beta(\kappa_s)=T_s$.
Then $\kappa_s$ is an $(s,s+1,s+2)$-core partition of maximum size.
Moreover, if $s$ is even, then $\kappa_s$ is the unique $(s,s+1,s+2)$-core partition of maximum size, which is self-conjugate.
 If $s $ is odd, then there is exactly another $(s,s+1,s+2)$-core partition of maximum size, which is the conjugate of $\kappa_{s}$.
\end{thm}

\pf
Let $\lambda$ be an $(s,s+1,s+2)$-core partition of maximum size.
We aim to show that $\lambda=\kappa_s$ if $s$ is even, and $\lambda=\kappa_s$ or $\kappa_s'$ if $s$ is odd.
From Lemma \ref{lem:triangleideal} we see that $\lambda=\lambda_{i,j}$ for some integers $i,j$.
To determine the values of $i,j$,
we consider the following two cases.

\noindent Case 1: $s$ is even.
 As a consequence of Lemma \ref{lem:triangle},
  we have $\lambda_{i,k}<\lambda_{i,\lfloor \frac{s-i+1}{2}\rfloor}$ for $0\leq k<\lfloor \frac{s-i+1}{2} \rfloor$.
Hence $j=\lfloor \frac{s-i+1}{2} \rfloor$, that is, $\lambda=\lambda_{i,\lfloor \frac{s-i+1}{2}\rfloor}$ for some $i$.
We claim that $i=1$.
Suppose to the contrary that $i>1$, that is, $i-1\geq 1$.
 By the definition of $\lambda_{i,j}$,  we find that $\lambda_{i,\lfloor \frac{s-i+1}{2}\rfloor}=\lambda_{i-1,0}$.
Since $s$ is even,
by Lemma \ref{lem:triangle}, we obtain that
\[
|\lambda|=
|\lambda_{i,\lfloor \frac{s-i+1}{2}\rfloor}|=|\lambda_{i-1,0}|<|\lambda_{i-1,\lfloor \frac{s-i+2}{2}\rfloor}|,
\]
contradicting the fact that $\lambda$ is of maximum size.
Hence we have $i=1$, and so $\lambda=\lambda_{1, \lfloor \frac{s}{2}\rfloor}$.

\noindent Case 2: $s$ is odd.
We claim that $i\leq 2$.
Suppose   that $i>2$.
By Lemma \ref{lem:triangle}, we have $j=\lfloor \frac{s-i+1}{2} \rfloor$, that is, $\lambda=\lambda_{i,\lfloor \frac{s-i+1}{2}\rfloor}$ for some $i$.
Since $i-1>1$ and $\lambda_{i,\lfloor \frac{s-i+1}{2}\rfloor}=\lambda_{i-1,0}$, using Lemma \ref{lem:triangle} we get
\[
|\lambda|=
|\lambda_{i,\lfloor \frac{s-i+1}{2}\rfloor}|=|\lambda_{i-1,0}|<|\lambda_{i-1,\lfloor \frac{s-i+2}{2}\rfloor}|,
\]
which contradicts the fact that $\lambda$ is of maximum size.
This proves the claim, namely, $i=1$ or $2$.
By Lemma \ref{lem:triangle}, we obtain  $\lambda=\lambda_{1,\lfloor \frac{s}{2}\rfloor}$, $\lambda_{1,0}$ or
$\lambda_{2,\lfloor \frac{s-1}{2}\rfloor}$.
By the definition of $\lambda_{i,j}$, we see that
$\lambda_{1,0}=\lambda_{2,\lfloor \frac{s-1}{2}\rfloor}$.
Thus, $\lambda=\lambda_{1,\lfloor \frac{s}{2}\rfloor}$ or
$\lambda_{1,0}$.
 Again, using Lemma \ref{lem:triangle}, we have
 $|\lambda_{1,\lfloor \frac{s}{2}\rfloor}|=|\lambda_{1,0}|$.
So we find that $\lambda_{1,\lfloor\frac{s}{2}\rfloor}$ and $\lambda_{1,0}$ are the only two $(s,s+1,s+2)$-core partitions of maximum size.

We now  conclude that
if $s$ is even, the partition $\lambda_{1,\lfloor\frac{s}{2}\rfloor}$ is the unique $(s,s+1,s+2)$-core partition of maximum size, and if $s$ is odd, there are two $(s,s+1,s+2)$-core partitions of maximum size, which are $\lambda_{1,\lfloor\frac{s}{2}\rfloor}$ and $\lambda_{1,0}$, respectively.

Notice that in both cases, $\lambda_{1,\lfloor\frac{s}{2}\rfloor}$ is an $(s,s+1,s+2)$-core partition of maximum size.
To prove that $\kappa_s$ is of maximum size, it suffices to show that $\kappa_s=\lambda_{1,\lfloor\frac{s}{2}\rfloor}$.
By the definitions of $\lambda_{i,j}$ and $\beta_{i,j}$, it can be verified that $\beta(\lambda_{1,\lfloor\frac{s}{2}\rfloor})
=\beta_{1,\lfloor\frac{s}{2}\rfloor}=T_s$.
Since $\beta(\kappa_s)=T_s$, we have $\beta(\kappa_s)=\beta(\lambda_{1,\lfloor\frac{s}{2}\rfloor})$.
This implies that $\kappa_s=\lambda_{1,\lfloor\frac{s}{2}\rfloor}$.
So we reach the conclusion that
 $\kappa_s$ is an $(s,s+1,s+2)$-core partition of maximum size.

It remains to show that
when $s$ is even,   $\kappa_s$ is self-conjugate, and when $s$ is odd,
 $\lambda_{1,0}=\kappa'_s$.

Clearly, the conjugate of an $(s,s+1,s+2)$-core partition
is still an $(s,s+1,s+2)$-core partition of the same size.
Since $\kappa_s$ is an $(s,s+1,s+2)$-core partition of maximum size,
 $\kappa_{s}'$ is also an
$(s,s+1,s+2)$-core partition of maximum size.

 When $s$ is even, since  $\kappa_s=\lambda_{1,\lfloor\frac{s}{2}\rfloor}$ is the unique $(s,s+1,s+2)$-core partition of maximum size,
 we have $\kappa_s'=\kappa_{s}$, that is,  $\kappa_s$ is self-conjugate.

When $s$ is odd, we have shown that $\kappa_{s}$ and $\lambda_{1,0}$
are the only two $(s,s+1,s+2)$-core partitions of maximum size.
To prove that $\lambda_{1,0}=\kappa_{s}'$,
it suffices to show that $\kappa'_s\not=\kappa_s$, that is,
$\kappa_s$ is not self-conjugate.
To this end, we aim to prove that the length of $\kappa_s$
is not equal to the largest part of $\kappa_s$.

Assume that $s=2m-1$ for some $m\geq 2$.
Note that the length of $\kappa_s$
equals $|\beta(\kappa_s)|=|T_s|$.
In view of  Theorem \ref{thm:poset},  we obtain that
\begin{equation}\label{size(Ts)}
|T_s|=\sum_{k=0}^{m-2}(2m-2-2k)=m^2-m.
\end{equation}
Thus the length of $\kappa_s$
equals $m^2-m$.
Since $\beta(\kappa_s)=T_s$,
from \eqref{betaset} and Theorem \ref{thm:poset},
it can be seen that the largest part of $\kappa_s$ equals
\[\left\lfloor\frac{s}{2}\right\rfloor s-1-(|T_s|-1)=(m-1)(2m-1)-(m^2-m)=m^2-2m+1.\]
Since $m\geq 2$, we have $m^2-m \neq m^2-2m+1$,
so that  the length of $\kappa_s$
is not equal to the largest part of $\kappa_s$.
This completes the proof.
\qed

 Theorem \ref{thm:maxcore} says that the partition $\kappa_s$   corresponding to the order ideal $T_s$ is of maximum size.
This leads to a proof of Conjecture \ref{conj2} which gives an explicit formula for the maximum size of an $(s,s+1,s+2)$-core partition.

\begin{cor}
Let $s$ be a positive integer. The maximum size $l(s)$ of an $(s,s+1,s+2)$-core partition equals
\begin{equation*}
l(s)=\left\{
       \begin{array}{ll}
         m{m+1\choose 3}, & \mbox{if } s=2m-1, \\[6pt]
         (m+1){m+1\choose 3}+{m+2\choose 3},  & \mbox{if } s=2m.
       \end{array}
     \right.
\end{equation*}
\end{cor}

\pf
It is easily checked that the corollary holds for $s\leq 2$.
We now assume that $s\geq 3$.
By Theorem \ref{thm:maxcore}, we know that the partition $\kappa_s$ such that $\beta(\kappa_s)=T_s$ is of maximum size.
Using   \eqref{equ:size}, we get
\begin{equation}\label{l(s)}
l(s)=|\kappa_s|=\sum_{h\in T_s}h-{|T_s| \choose 2}.
\end{equation}

If $s$ is odd, that is, $s=2m-1$ for some $m\geq 2$,  by Theorem \ref{thm:poset}, we find that
\begin{align}\label{sum(Ts)}
\sum_{h\in T_s}h
&=\sum_{k=0}^{m-2}\ \sum_{i=1}^{2m-2k-2}\big(k(2m+1)+i\big)\notag\\[3pt]
&=\sum_{k=0}^{m-2}(-4mk^2+4m^2k-6mk+k+2m^2-3m+1)\notag\\[3pt]
&=\frac{2}{3}m^4-m^3-\frac{1}{6}m^2+\frac{1}{2}m.
\end{align}
Substituting \eqref{size(Ts)} and \eqref{sum(Ts)} into \eqref{l(s)}, we obtain that
\[
l(s)=m{m+1\choose 3}.
\]

If $s$ is even, that is, $s=2m$ for some $m\geq 2$,
by Theorem \ref{thm:poset}, we obtain that
\begin{align}
\sum_{h\in T_s}h
&=\sum_{k=0}^{m-1}\ \sum_{i=1}^{2m-2k-1}\big(k(2m+2)+i\big)\notag\\[3pt]
&=\sum_{k=0}^{m-1}(-4mk^2-2k^2+4m^2k-2mk-k+2m^2-m)\notag\\[3pt]
&=\frac{2}{3}m^4+\frac{1}{3}m^3-\frac{1}{6}m^2+\frac{1}{6}m.\label{sum(T2m)}
\end{align}
Again, by Theorem \ref{thm:poset}, we get
\begin{equation}\label{size(T2m)}
|T_s|=\sum_{k=0}^{m-1}(2m-2k-1)=m^2.
\end{equation}
Substituting  \eqref{sum(T2m)} and \eqref{size(T2m)} into \eqref{l(s)} gives
\[
l(s)=(m+1){m+1\choose 3}+{m+2\choose 3}.
\]
This completes the proof.
\qed

\section{Proof of Conjecture \ref{conj3}}\label{sec:4}

In this section, we shall give a proof of Conjecture \ref{conj3} on the total sum $h(s)$ of sizes of   $(s,s+1,s+2)$-core partitions.
By the correspondence between $(s,s+1,s+2)$-core partitions and
the order ideals of $T_s$,
we can express $h(s)$
in terms of the sums of elements of order ideals of $T_s$. Then
we obtain an explicit formula for the generating function of $h(s)$, which leads to a proof of Conjecture \ref{conj3}.

By Theorem \ref{thm:general} and relation \eqref{equ:size}, we have
\begin{equation}\label{equ:h}
h(s)=\sum_{I\in J(T_s)}\left(\sum_{a\in I}a-{|I|\choose 2}\right).
\end{equation}
Let $\rho_s$ denote the rank function of the poset $T_s$.
By Theorem \ref{thm:poset}, we see that
$ \rho_s(a)=k$ for $a\in B_k=[1+k(s+2),(k+1)s-1]$.
In order to derive the generating function of $h(s)$, we need
the following two functions
\begin{align*}
f(s)&=\sum_{I\in J(T_s)}|I|,\\[6pt]
g(s)&=\sum_{I\in J(T_s)}\sum_{a\in I}\rho_s(a).
\end{align*}
Let $F(x)$, $G(x)$ and $H(x)$ be the ordinary generating functions of the numbers $f(s)$, $g(s)$ and $h(s)$,
that is,
\begin{align*}
    F(x)&=\sum_{s\geq 0}f(s)x^{s},\\[3pt]
    G(x)&=\sum_{s\geq 0}g(s)x^{s},\\[3pt]
    H(x)&=\sum_{s\geq 0}h(s)x^{s}.
\end{align*}

The following lemma gives  recurrence relations for $f(s)$, $g(s)$ and $h(s)$,  which lead to the generating functions $F(x)$, $G(x)$ and $H(x)$.
In fact, we shall use the generating functions
 $F(x)$ and $G(x)$ to compute $H(x)$.

\begin{lem}\label{lem:recur}
For $s\geq 2$, we have
\begin{align}
f(s)=
&f(s-1)+\sum_{i=2}^s\big(2M_{s-i}f(i-2)+(i-1)M_{s-i}M_{i-2}\big),\label{rec:fs}\\[3pt]
g(s)=
&g(s-1)+\sum_{i=2}^sM_{s-i}\big(2g(i-2)+f(i-2)\big),\label{rec:gs}\\[3pt]
h(s)=&h(s-1)+f(s-1)+g(s-1)\notag\\[3pt]
    &\quad +\sum_{i=2}^s\big(2M_{s-i}h(i-2)+(s+4-i)M_{s-i}f(i-2)\notag\\[3pt]
      &\quad \qquad +2(s-i+2)M_{s-i}g(i-2)+(i-1)M_{s-i}M_{i-2}
      -f(i-2)f(s-i)\big).\label{rec:hs}
\end{align}
\end{lem}

\pf
We shall only give a proof of \eqref{rec:hs}.
Relations \eqref{rec:fs} and \eqref{rec:gs} can be verified in the same manner.

Let
\[
h_i(s)=\sum_{I\in J_i(T_s)}\left(\sum_{a\in I}a-{|I|\choose 2}\right).
\]
Since
\[
J(T_s)=\bigcup_{i=1}^sJ_i(T_s),
\]
in view of \eqref{equ:h}, we have
\[
h(s)=\sum_{i=1}^sh_i(s).
\]

To compute $h_i(s)$, we recall the decomposition of an order ideal of $T_s$ as given in Section 2.
For an order ideal $I\in J_i(T_s)$, we can express $I$ as
\begin{equation}\label{equ:decompsition}
I = \{1,2,\ldots,i-1\}\cup I'\cup I'',
\end{equation}
where $I'$ is isomorphic to an order ideal $I_1$ of $T_{i-2}$ and $I''$ is isomorphic to an order ideal $I_2$ of $T_{s-i}$.
Here we set $T_{-1}$ to be the empty set.
Conversely,  an order ideal $I_1$ of $T_{i-2}$ and an order ideal $I_2$ of $T_{s-i}$ uniquely determine an order ideal $I\in J_i(T_s)$.
To be more specific, we have
\[
I'=\{a+s+2+(s+2-i)\rho_{i-2}(a)\mid a\in I_1\}
\]
and
\[
I''=\{a + i + i\rho_{s-i}(a)\mid a\in I_2\}.
\]
The above decomposition implies that for $I\in J_i(T_s)$,
\begin{align*}
\sum_{a\in I}a
&=\binom{i}{2}+\sum_{a\in I'}a+ \sum_{a\in I''}a\\[6pt]
&=\binom{i}{2}+\sum_{a\in I_1} \big(a+s+2+(s+2-i)\rho_{i-2}(a)\big)+ \sum_{a\in I_2}\big(a + i + i\rho_{s-i}(a)\big).
\end{align*}
From the proof of Conjecture \ref{conj1}, we see that
$M_s$ equals the number of order ideals of $T_s$.
Let \[ p(s)=\sum_{I\in J(T_s)}\sum_{a\in I}a.\]
Using the decomposition \eqref{equ:decompsition} of an
order ideal $I\in J_{i}(T_s)$,  $h_{i}(s)$ can be computed as follows.
For $i=1$, we have
\begin{align}
h_1(s)&=\sum_{I\in J_1(T_s)}\left(\sum_{a\in I}a-{|I|\choose 2}\right)\notag\\[6pt]
&=\sum_{I_2\in J(T_{s-1})}\left(\sum_{a\in  I_2}\big(a + 1 + \rho_{s-1}(a)\big)-
\binom{|I_2|}{2}\right)\notag\\[6pt]
&=h(s-1)+f(s-1)+g(s-1).
\end{align}
For $2\leq i\leq s$, we find that
\begin{align}\label{equ:gis}
h_i(s)=&\sum_{I\in J_i(T_s)}\left(\sum_{a\in I}a-\binom{ |I|}{2}\right)\notag\\[6pt]
=&\sum_{\substack{I_1\in J(T_{i-2})\\[1pt]I_2\in J(T_{s-i})}}
\left(\binom{i}{2}+\sum_{a\in I_1} \big(a+s+2+(s+2-i)\rho_{i-2}(a)\big)\right.\notag\\[6pt]
 &\qquad\left.+ \sum_{a\in I_2}\big(a + i + i\rho_{s-i}(a)\big)-\binom{ |I_1|+|I_2|
  +i-1}{2}\right)\notag\\[6pt]
=&\sum_{I_1\in J(T_{i-2})}M_{s-i}\left(\sum_{a\in
  I_1}\big(a+s+2+(s+2-i)\rho_{i-2}(a)\big)+\binom{i}{2}\right) \notag\\[6pt]
  &\quad+ \sum_{I_2\in
  J(T_{s-i})}M_{i-2}\sum_{a\in I_2}\big(a + i + i\rho_{s-i}(a)\big)\notag\\[6pt]
  &\quad-\sum_{\substack{I_1\in J(T_{i-2})\\[1pt]I_2\in J(T_{s-i})}}\binom{ |I_1|+|I_2|
  +i-1}{2} \notag\\[6pt]
 =&M_{s-i}\big(p(i-2)+(s+2)f(i-2)+
(s-i+2)g(i-2)\big)\notag\\[6pt]
&\quad+M_{s-i}M_{i-2}\binom{i}{2}
 +M_{i-2}\big(p(s-i)+if(s-i)+ig(s-i)\big)\notag\\[6pt]
 &\quad- \sum_{\substack{I_1\in J(T_{i-2})\\[1pt]I_2\in J(T_{s-i})}}\binom{|I_1|+|I_2|
  +i-1}{2}.
\end{align}
Since
\begin{align*}
\sum_{I_1\in J(T_{i-2})}|I_1|=f(i-2), \\[3pt]
\sum_{I_2\in J(T_{s-i})}|I_2|=f(s-i),
\end{align*}
we have
\begin{align}\label{equ:binomial}
&\sum_{\substack{I_1\in J(T_{i-2})\\[1pt]
                   I_2\in J(T_{s-i})}}
     \binom{|I_1|+|I_2|+i-1}{2}\notag\\[6pt]
=& \sum_{\substack{I_1\in J(T_{i-2})\\[1pt]
                   I_2\in J(T_{s-i})}}
   \left(\binom{|I_1|}{2}+\binom{|I_2|}{2}+(i-1)|I_1| +(i-1)|I_2|+|I_1||I_2|+\binom{i-1}{2}\right)\notag\\[6pt]
=& \sum_{I_1\in J(T_{i-2})}
     M_{s-i}\left(\binom{|I_1|}{2}+(i-1)|I_1|\right)
   +\sum_{I_2\in J(T_{s-i})}
     M_{i-2}\left(\binom{|I_2|}{2}+(i-1)|I_2|\right)\notag\\[6pt]
&\quad +\bigg(\sum_{I_1\in J(T_{i-2})}|I_1|\bigg)
        \bigg(\sum_{I_2\in J(T_{s-i})}|I_2|\bigg)
       +M_{s-i}M_{i-2}\binom{i-1}{2}\notag\\[6pt]
=
& M_{s-i}\sum_{I_1\in J(T_{i-2})}\binom{|I_1|}{2}+(i-1)M_{s-i}f(i-2)
   +M_{i-2}\sum_{I_2\in J(T_{s-i})}\binom{|I_2|}{2}\notag\\[6pt]
&\quad +(i-1)M_{i-2}f(s-i)+f(i-2)f(s-i)+M_{s-i}M_{i-2}\binom{i-1}{2}.
\end{align}
Note that
\[h(s)=p(s)-\sum_{I\in  J(T_{s})}\binom{|I|}{2},\]
Substituting (\ref{equ:binomial}) into (\ref{equ:gis}),
 we obtain that
\begin{align*}
h_i(s)=
& M_{s-i}\big(p(i-2)+(s+2)f(i-2)+
     (s-i+2)g(i-2)\big)\\[6pt]
& \quad +M_{s-i}M_{i-2}\binom{i}{2}
        +M_{i-2}\big(p(s-i)+if(s-i)+ig(s-i)\big)\\[6pt]
& \quad -M_{s-i}\sum_{I_1\in J(T_{i-2})}\binom{|I_1|}{2}
        -(i-1)M_{s-i}f(i-2)\\[6pt]
& \quad -M_{i-2}\sum_{I_2\in J(T_{s-i})}\binom{|I_2|}{2}
        -(i-1)M_{i-2}f(s-i)\\[6pt]
& \quad -f(i-2)f(s-i)-M_{s-i}M_{i-2}\binom{i-1}{2}\\[6pt]
=
& M_{s-i}h(i-2)+(s+3-i)M_{s-i}f(i-2)\\[6pt]
& \quad +(s-i+2)M_{s-i}g(i-2)+(i-1)M_{s-i}M_{i-2}\\[6pt]
& \quad+M_{i-2}h(s-i)+M_{i-2}f(s-i)\\[6pt]
& \quad +iM_{i-2}g(s-i)-f(i-2)f(s-i).
\end{align*}
Summing over $i$, we deduce that
\begin{align*}
h(s)=&\sum_{i=1}^{s} h_i(s)\notag\\[3pt]
    =&h(s-1)+f(s-1)+g(s-1)\notag\\[3pt]
     &\quad +\sum_{i=2}^{s} \big(M_{s-i}h(i-2)+(s+3-i)M_{s-i}f(i-2)\notag\\[3pt]
     &\qquad +(s-i+2)M_{s-i}g(i-2)+(i-1)M_{s-i}M_{i-2}\notag\\[3pt]
    &\qquad +M_{i-2}h(s-i)+M_{i-2}f(s-i)\notag\\[3pt]
    &\qquad +iM_{i-2}g(s-i)-f(i-2)f(s-i)\big)\notag\\[3pt]
    =&h(s-1)+f(s-1)+g(s-1)\notag\\[3pt]
    &\quad +\sum_{i=2}^s\big(2M_{s-i}h(i-2)+(s+4-i)M_{s-i}f(i-2)\notag\\[3pt]
    &\qquad +2(s-i+2)M_{s-i}g(i-2)\notag\\[3pt]
      &\qquad +(i-1)M_{s-i}M_{i-2}-f(i-2)f(s-i)\big).
\end{align*}
This proves \eqref{rec:hs}.
\qed

From the recurrence relations   in Lemma \ref{lem:recur}, we get the following explicit formula for the generating function $H(x)$.

\begin{thm}\label{thm:gen_h(x)}
We have
\begin{align}\label{gef:G(X)}
H(x)={\frac {{x}^{2}}{ ( 1-2x-3{x}^{2} ) ^{5/2}}}.
\end{align}
\end{thm}

\pf
Let $M(x)$ be the ordinary generating function of the
Motzkin numbers, that is,
\[
M(x)=\sum_{s\geq 0}M_sx^s.
\]
It follows from  \eqref{rec:hs} that
\begin{align}\label{gef1:G(x)}
H(x)=&xH(x)+xF(x)+xG(x)+2x^2M(x)H(x)+4x^2M(x)F(x)\notag\\[3pt]
     &    \quad+x^3M'(x)F(x)
          +4x^2 M(x)G(x)+2x^3M'(x)G(x)\notag\\[3pt]
     &    \quad+x^3M'(x)M(x)+x^2M^2(x)-x^2F^2(x).
\end{align}
It is known that
\begin{align}\label{gef:T(x)}
    M(x)=\frac{1-x-\sqrt{1-2x-3x^2}}{2x^2}.
\end{align}
To derive a formula for $H(x)$, we proceed to compute
$F(x)$ and $G(x)$.

From recurrence relation (\ref{rec:fs}) we get
\begin{align}\label{rec:H(x)}
F(x)=xF(x)+2x^2M(x)F(x)+x^3\left(M'(x)
+\frac{M(x)}{x}\right)M(x).
\end{align}
 Substituting (\ref{gef:T(x)}) into (\ref{rec:H(x)}), we obtain
\begin{align}\label{gef:H(x)}
F(x)={\frac { \left( -1+x+\sqrt {1-2x-3{x}^{2}} \right) ^{2}}{4{x}
^{2} ( 1-2x-3{x}^{2} ) }}.
\end{align}

Similarly, from recurrence relation \eqref{rec:gs} we deduce that
\begin{align*}
G(x)=xG(x)+2x^2M(x)G(x)+x^2M(x)F(x),
\end{align*}
which implies that
\begin{align}\label{gef1:F(x)}
G(x)=\frac{x^2M(x)F(x)}{1-x-2x^2M(x)}.
\end{align}
Substituting (\ref{gef:T(x)}) and (\ref{gef:H(x)}) into (\ref{gef1:F(x)}), we get
\begin{align}\label{gef:F(x)}
G(x)=-{\frac { \left( -1+x+\sqrt {1-2x-3{x}^{2}} \right) ^{3}}{8{x}
^{2} ( 1-2x-3{x}^{2} ) ^{3/2}}}.
\end{align}
Based on the formulas for $F(x)$ and $G(x)$,
the formula \eqref{gef:G(X)} for $H(x)$ immediately follows from (\ref{gef1:G(x)}). This completes the proof.
\qed

Theorem \ref{thm:gen_h(x)} leads to   an explicit
formula for $h(s)$, which confirms Conjecture \ref{conj3}.

\begin{cor}
Let $s$ be a positive integer. The sum of the sizes of all the $(s,s+1,s+2)$-core partitions is
\begin{equation}\label{equ:h(s)}
h(s)=\sum _{j=0}^{s-2} {j+3\choose 3}\sum _{i=0}^{\lfloor j/2 \rfloor}
{j\choose 2i}C_i.
\end{equation}
\end{cor}

\pf
Exchanging the order of summations on the
right hand side of \eqref{equ:h(s)}, we get
\begin{equation}\label{exp:h(s)}
h(s)=\sum_{i\geq 0}\frac{(2i+3)!}{6\, i!(i+1)!}{s+2\choose 2i+4}.
\end{equation}
On the other hand,  based on the expression for $H(x)$ as in
 (\ref{gef:G(X)}),
 it can be verified that  $H(x)$
satisfies the following differential equation
\begin{equation*}
(-x+2x^2+3x^3)H'(x)+(2+x+9x^2)H(x)=0,
\end{equation*}
which implies that for $s\geq 3$,
\begin{equation}\label{equ:rec_h(s)}
(2-s)h(s)+(2s-1)h(s-1)+(3s+3)h(s-2)=0.
\end{equation}
Using the Zeilberger algorithm, see \cite{PWZ},
we find that the sum in
\eqref{exp:h(s)} also
satisfies the same recurrence relation \eqref{equ:rec_h(s)}
 as $h(s)$. Taking the initial values into consideration,
 we arrive at \eqref{exp:h(s)}. This completes the proof.
\qed

After the completion of this work, we noticed that
part of the results in this paper  were independently
obtained by Amdeberhan and Leven \cite{A-L}.

\vspace{0.5cm}
 \noindent{\bf Acknowledgments.}
 This work was supported by the 973 Project, the PCSIRT Project of the Ministry of Education and the National Science Foundation of China.

\end{document}